\documentclass[11pt]{article}

\usepackage{amsmath,latexsym,amscd,amsthm,amssymb}
\usepackage{times} 
\usepackage{mathptm}

\newtheorem{Thm}{Theorem}[section]
\newtheorem{Q}[Thm]{Question}

 \newtheorem{Lma}[Thm]{Lemma}
 \newtheorem{Def}[Thm]{Definition}
 \newtheorem{Cor}[Thm]{Corollary}
 \newtheorem{Rem}[Thm]{Remark}
  
 \newtheorem{Prop}[Thm]{Proposition}

 \newcommand{\Hom}{{\rm Hom}}
 \newcommand{\height}{{\rm ht}}
 \newcommand{\depth}{{\rm depth}}
 
 \newcommand{\Ann}{{\rm Ann}}

\newcommand{\into}{\hookrightarrow}
\newcommand{\Dim}{{\rm dim}}

\newcommand{\tc}{0^{\ast}}

\newcommand{\er}{0^{\ast fg}_{E_{R}}}
 \newcommand{\hr}{0^{\ast} _{H_{m}^{d}(R)}}
\newcommand{\fghr}{0^{\ast fg} _{H_{m}^{d}(R)}} 
\newcommand{\tchi}{0^{\ast} _{H_{m}^{d}(I)}}
 \newcommand{\fgm}{0^{\ast fg}_M}
\newcommand{\Spec}{{\rm Spec}}

\newcommand{\taup}{\tau_{par}}
\newcommand{\tauo}{\tau_{par}(\omega_R)}

\begin{document}
 
\title{\bf Strong Test Modules and Multiplier Ideals}
\author{Florian Enescu \thanks{Department of Mathematics, University of Utah, Salt Lake City, Utah, 84112, USA, and the Institute of Mathematics of the Romanian Academy, Bucharest, Romania. MSC: 13A35}}
\date{}

\maketitle

{\bf Abstract.} We introduce the notion of strong test module and show that a large number of such modules appear in the tight closure theory of complete domains: the test ideal (this has already been known), the parameter test module, and the module of relative test elements. They also appear as certain multiplier ideals, a concept of interest in algebraic geometry.

\section {\bf Introduction}

The purpose of this note is to address a few issues related to the tight closure of ideals in rings of characteristic $p > 0$. The study regards the concept of strong test ideals introduced by C.~Huneke (\cite{hu}). A.~Vraciu (\cite{vr}), N.~Hara and K.~E.~Smith (\cite{ha-sm}) have also investigated it. In this note, we extend the notion of strong test ideals to modules and generalize some of the known results. The perspective that we offer in our study also leads to some algebraic properties that are shared by the multiplier ideals. Along the way, we provide a natural interpretation for the modules of relative test elements, a concept introduced by A.~K.~Singh (\cite{si}).

Throughout these notes, $(R,m,K)$ is a local ring of characteristic $p>0$. 
We start by recalling the notions of tight closure for modules, test ideal and strong test ideal. 

In positive characteristic $p$, one can define the Frobenius homomorphism $F : R \to R$, $F(r)=r^p$. For every $e$, the iterated Frobenius map $F^e : R \to R$ sends $r$ to $r^{p^e}$ and enables $R$ with a new $R$-algebra structure on the right, denoted by $R^e$ (on the left, $R^e = R$). It also defines a functor that sends an $R$-module $M$ to $F^e(M):= R^e \otimes_R M$. 

Let $M$ an $R$-module and $N \subset M$ a submodule in $M$. The {\it tight closure of $N$ in $M$}, denoted by $N^\ast_M$, is defined as follows: $m \in N^\ast_M$ if there is $c \in R^o: =R - \cup_{P \in Min(R)} P$ such that $c \otimes m \in N^{[q]}_M:={\rm Im}(F^e(N) \to F^e(M))$, for all $e$ sufficiently large. The element $c \otimes m$ belongs to $F^e(M)$ and it is ocasionally denoted by $cm^{[q]}$, where $q=p^e$. Whenever $M=R$ and $N=I$ an ideal of $R$, we simply obtain the tight closure of $I$ in $R$, denoted by $I^\ast$. In the definition of tight closure, there is no assumption that $M$ is finitely generated over $R$. In fact, there is another
notion of tight closure called {\it the finitistic tight closure}. If $N$ is a submodule of $M$, the finitistic tight closure of $N$ in $M$ is denoted by $N^{\ast fg}_M$ and equals $\cup _{M'} (N \cap M')^\ast _{M'}$, where the union runs over all finitely generated submodules $M'$ of $M$. It is easy to see that when $M$ is finitely generated over $R$, $N^\ast_M = N^{\ast fg}_M$. The case when $M$ is Artinian and $N=0$ is especially important. It is conjectured that in this case the two notions of tight closure coincide (the conjecture has been proven in a few cases, \cite{ly-sm}). An important special case of the conjectured is that when $M = E_R(k)$.
For basic tight closure facts we refer the reader to~\cite{ho-hu}, where the theory is presented in detail. For details on the above stated conjecture on Artinian modules and related issues we refer the reader to~\cite{ly-sm}.

\begin{Def}
{\rm Let $R$ be as above. The} test
ideal {\rm $\tau(R)$ is defined as $\cap_{M} \Ann_{R} (0^\ast _{M})$,
where $M$ runs through all finitely generated $R$-modules. An element of $\tau(R) \cap R^o$ is called } a test element.
\end{Def}

We now list the main facts about the test ideal.

\begin{Prop}
\label{prop-test}
Let $R$ be a Noetherian ring of characteristic $p$. Then

(a) $c \in \tau(R)$ if and only if whenever $N \subset M$ and $x \in
N^\ast _{M}$, then $cx^q \in N^{[q]}_{M}$ for all $q$. In fact $c \in R^o$ 
is a test element if and only if $cI^\ast \subset I$, for every ideal $I$ in $R$.

(b) $\tau(R) = \Ann_{R} (\er)$, where $E_{R} = \oplus_{m} E(R/m)$, and
$m$ runs through all maximal ideals of $R$.

(c) If  $(R,m,K)$ is local and $\{ I_{t} \}$ is a sequence of
$m$-primary irreducible ideals cofinal with the powers of $m$, then
$\tau (R) = \cap_{t} I_{t} : I_{t}^\ast$. Such as sequence exists if
and only if $R$ is approximately Gorenstein.  If $R$ is local and
complete, then $\Ann_{E_{R}}{(\tau(R))} = \er$.

\end{Prop}

The notion of the strong test ideal has been introduced by C.~Huneke as follows:

\begin{Def}
{\rm An ideal $T$ of $R$ such that $T I^\ast =T I$ for every ideal
$I$ is called} a strong test ideal.
\end{Def}

The motivation for this definition can be explained as follows. Since $I^\ast \subset \overline{I}$, then for every element $x \in I^\ast$ there is an integral dependence equation over $I$ that is satisfied by $x$. If $T$ is a strong test ideal and $R$ is domain, then an application of the determinant trick (see Theorem 2.1 in~\cite{hu}) shows that the minimal number of generators of $T$ provides an upper bound for the degree of such an integral dependence equation. What is significant here is that this bound is uniform for every ideal $I$ and every $x \in I^\ast$, as it depends only on the ideal $T$. In the general case, the existence of such bound can be reduced to the domain case. Finding more than one strong test ideal is important in practice as we are not aware of any result that indicate which one has fewest number of generators. In addition to this, there is another aspect of the definition. Each time a strong test ideal is exhibited, its defining property gives a uniform special feature of tight closure of ideals. In some cases, this can be useful in applications (Theorem 5.2 in~\cite{ab-en}).

\section{\bf Test Modules and Multiplier Ideals}

In this section we define the notion of strong test module and give examples of such modules.

A.~Vraciu has proven an important property of the test ideal in complete rings.

\begin{Thm}[Vraciu]
Let $(R,m,K)$ be a complete local ring. Then $\tau(R)$ is a strong
test ideal.
\end{Thm}

We will provide a natural generalization of this theorem and link it to multiplier ideals via the results of N.~Hara (\cite{ha}).

Throughout this section $R$ is assumed complete. Also, for simplicity, $R$ will  be assumed domain in some results. As illustrated in ~\cite{hu}, the issues related to strong test ideals can be reduced to the domain case in many instances, by reduction modulo each minimal prime.

Let $M$ be an Artinian module over $R$ and denote by $M^{\vee}$ the
Matlis dual $\Hom_R (M,E_R(K))$. Clearly, the duality induces a natural pairing:

$$M^\vee \times M \to E_R(K). \hspace{2cm} (\ast) $$

Using this pairing, let us define $ \tau _M: = \{ n \in M^\vee : n \tc _M =0 \}.$ A similar definition can be formulated using $\fgm$ instead of 
$\tc _M$. A module of this type will be denoted $\tau^{fg}_M$. A module of the form $\tau_M$ or $\tau^{fg}_M$ will be called {\it a test module}.
 It is likely that $\tau_M =\tau^{fg}_{M}$ in general, as it has been conjectured that $\fgm = \tc _M$ for an Artinian module $M$ (see Theorem 8.12 in~\cite{ly-sm} where the conjecture is proven for local rings with isolated singularities; the important case $M=E_R(k)$ it also settled there for Cohen-Macaulay local rings, Gorenstein on their punctured spectrum, Theorem 8.8,~\cite{ly-sm}).

We would like to state the following useful fact from Matlis duality theory (see Lemma 2.1 in~\cite{sm2} or Lemma 3.3 in~\cite{ha}).

\begin{Lma}
\label{ann}
Let $F$ be a finitely generated module over a local ring $(R,m)$. Denote its Matlis dual by $F^\vee = \Hom_R(F, E_R)$. Let $G$ (resp. $L$) be any submodule of $F$ (resp. $F^\vee$). If $L=\Ann_{F^\vee} G$, then $G=\Ann_{F} L$. If $(R,m)$ is complete, then the reverse is also true. 
\end{Lma}

\begin{Thm}
\label{strong}
Let $R$ be a local complete ring and $M$ an Artinian module over $R$  defining the test modules $\tau _M$ and $\tau^{fg}_M$. Then $I^\ast \tau _M = I \tau _M$ and $I^\ast \tau^{fg}_M = I \tau^{fg}_M $ 
for every ideal $I$ of $R$.
\end{Thm}

\begin{proof}

To prove the claimed equality it is enough to show that the two modules in the statement of the theorem have the
same annihilator in $M$ (here, we need $R$ be complete for local duality).

We will start with $\tau_M$.

$\Ann_M (I^\ast \tau) = \{ m \in M : I^\ast \tau m = 0 \} = \{ m \in M: 
I^\ast \cdot m \subset \Ann_M (\tau) \}$.  Local duality gives that $\Ann_M
(\tau) = \tc _M$. Therefore, $\Ann_M (I^\ast \tau) = \{ m: m I^\ast
\subset \tc _M \} =  (\tc _M : I^\ast). $ Similarly, $\Ann_M (I \tau) =
(\tc_M : I)$. So, we need to show that
$$(\tc _M : I ^\ast) = (\tc _M: I).$$

The inclusion $(\tc _M : I ^\ast) \subset (\tc _M: I)$ is evident, so we will concentrate on the reverse inclusion:

Let us take $z \in I^\ast$ and $m \in (\tc _M
: I)$. We need to show that $m \in (\tc _M : I ^\ast)$ and hence it suffices to show that
$zm \in \tc _M$.

Now, $z \in I^\ast$ so there is $d \in R^o$ such that $dz^q = \sum
a_{iq} x_{i}^q$, where $I = (x_1,...,x_k)$. Since $m \in M$, we get an element $m^q : = 1 \otimes m \in F^e(M)$ for every $q=p^e$. With this notation, 
 $dz^q \cdot m^q = \sum a_{iq} x_{i}^q \otimes m = \sum a_{iq}\otimes x_{i} m$. Since each $x_i \in I$ and $m \in (\tc _M
: I)$, one has that $x_i m \in \tc _M$. Take $c \in R^o $ an element which works the tight closure equations given by $x_i m \in \tc _M$ for {\it every} $i=1,...,k$. Then $c \otimes x_im = 0$ in $F^e(M)$, for all $i$. So, $cd(zm)^q = \sum a_{iq}c \otimes x_im = 0$ so $zm \in \tc_M$.

In the case of $\tau^{fg}_M$ the reasoning is similar. Keeping the notations as above, we have now that $x_im \in \fgm$ for every $i=1,...,k$. This means that there is a finitely generated submodule $M'$ of $M$ such that $x_im \in 0^\ast_{M'}$ for all $i$.
Now, we can use $c$ a test element of $R$ to show, as above, that $zm \in \fgm$.

\end{proof}

\begin{Def}
\label{strongfg}
 {\rm  Let $T$ be an $R$-module. The property that $I T = I^\ast T$, for every ideal $I$ in $R$, will be called} the strong test module property. {\rm A faithful module $T$ with the strong test module property is called a} strong test module. 
\end{Def}

\begin{Rem}{\rm The minimal number of generatos of a strong test module $T$ provides a uniform bound (depending only on the module $T$) on the degree of the equation of integral dependence that an element $x \in I^\ast$ satisfies over $I$, for every $I$ and every such $x$. This can be obtained by a straightforward generalization of the argument given in Theorem 2.1~\cite{hu}.}

\end{Rem}
Two special cases of the Theorem stand out. The first part of the next Corollary recovers Vraciu's result. The second part refers to the {\it parameter test module} $\tauo$, a notion introduced by Karen E.~Smith (~\cite{sm2}). Let us recall that $\tauo = \Ann_{\omega _R}(\hr)$ for an excellent local Cohen-Macaulay ring.

\begin{Cor}
Let $R$ as above.

1) The test ideal $\tau (R)$ is a strong test ideal.

2) Assume that $R$ is also  Cohen-Macaulay and domain. The parameter test module $\tauo$ has the strong test module property.
\end{Cor}

\begin{proof}
The Theorem~\ref{strong} applies in both cases. Also, $\tauo$ is a faithful module as a submodule of $\omega _R$ which is torsion-free in our case.

\end{proof}

Now, we take a look at {\it the module of relative test elements} for a finite extension of reduced $F$-finite local rings $(A,m_o, k) \to (R,m,K)$. The concept was introduced by Anurag K.~Singh in ~\cite{si} and is defined as $T(R,A):={\rm Ann}_{M^\vee} (\fgm)$ for $M = E_A(k) \otimes_A R$ under the duality $(\ast)$. Therefore, the module of relative test elements for $A \to R$ is a particular type of test module for $R$. It is worth noticing that $\Hom_R (E_A(k) \otimes_A R, E_R(K)) \simeq \Hom_A(R,A)$ and that $T(R,A) \subset \Hom_A(R, A).$ 

\begin{Cor}
Let $(A,m_o, k) \to (R,m,K)$  be a finite extension of $F$-finite local rings. Assume that $R$ is complete and domain. Then the module of relative test elements $T(R,A)$ is a strong test module.
\end{Cor}

\begin{proof}
The strong test module property of $T(R,A)$ follows immediately from Theorem~\ref{strong}. To show that $T(R,A)$ is faithful we would like to remark that $\Hom_A(R,A)$ is torsion-free, if $R$ is domain. It is enough to show that $rR \cap A \neq 0$, for every $r \in R$ (this implies that each $f \in \Hom_A(R,A)$ is injective). Since $R$ is module finite over $A$ it follows that $r$ in integral over $A$. Hence, there is an equation of minimal degree of the form $r^n + a_{n-1}r^{n-1}+...+a_1r + a_o =0$. 
But $R$ is a domain and $n$ has been chosen minimal, so $0 \neq a_o \in rR \cap A$.

\end{proof}

In some cases, the multiplier ideals can occur as test modules. In what follows we explain this assertion. We make use of the results of N.~Hara (\cite{ha})
who proved, in particular, that the test ideal is a certain multiplier ideal (this has also been proved independently by K.~E.~Smith,~\cite{sm}).

 First we need to describe the setup which is needed to state Hara's results. It involves reduction to positive characteristic from characteristic zero. The set up will be described without any proofs. All the assertions are addressed in detail in~\cite{ha} and  the reader should consult his paper (sections 4.6, 5.1 and 5.7).  

Let $R$ be a finitely generated algebra over a field of characteristic zero and let $I$ be a divisorial ideal such that $I^{(n)} \simeq R$ for some $n \in \mathbf{N}$. Consider $D$ a $\mathbf{Q}$-Cartier Weil divisor on $\Spec(R)$ such that $H^0(\Spec(R),\mathcal{O}(D)) = I$. The round-up of $D$ is denoted by $\ulcorner D \urcorner$ and the round-down by $\llcorner D \lrcorner$. Let $f: X \to \Spec (R)$ be a desingularization with exceptional divisor with simple normal crossing. 

Fix an isomorphism $I^{(n)} \simeq R$ and define two cyclic coverings:

$S = \bigoplus _{i=0}^{n-1} I^{(i)}$ and $ Y = {\it Spec}_{X}( \bigoplus_{i=0}^{n-1} \mathcal{O}_X(\llcorner if^*D \lrcorner))$. Also  let $h: \tilde{Y} \to Y$
be a resolution of singularities of $Y$. It is known that $Y$ has only rational singularities.

We have the following commutative diagram:

\begin{displaymath}
\begin{array}{ccccc}
\tilde{Y} & \to  & Y & \to & \Spec(S) \\

          &       & \downarrow &    & \downarrow \\

          &       &  X         &  \to & \Spec(R) \\

\end{array}
\end{displaymath}

We reduce all the data to characteristic $p \gg 0$ (and also localize at a prime ideal of the new algebra whenever we refer to the local case). As part of the set-up, we can assume that the above diagram is defined over a perfect field $\mathcal{K}$ of characteristic $p > 0$
which does not divide $n$. We will keep the notations unchanged.

If $(R,m)$ is local, $S$ is semilocal. We will denote $n_i$, $i=1,...,s$, the maximal ideals of $S$ and $n = \cap _i n_i$. $Z$ denotes the fiber of $X \xrightarrow{f} \Spec(R)$ and $Z_i$ denotes the fiber of $Y \xrightarrow{g} \Spec(S)$ over $n_i$, and $Z' = \cup _i Z_i$.

It can be seen that  $H^0(Y, \omega _Y) = \oplus _{i=0} ^{n-1} H^0(X,\omega _X ( \ulcorner -i f^\ast D \urcorner))$ and $\omega _S = \oplus_{i=0}^{n-1} \Hom(I^{(i)}, \omega _R)$ (these formulae differ by a sign from those in ~\cite{ha}, as this is what comes out from the direct application of the adjunction formula).

\begin{Thm}(Hara)
\label{hara}
Assume that $(R,m)$ is local, normal and of dimension $d$ which is obtained from characteristic zero by reduction to characterstic $p>0$ as above.
Then we have that 

$$ \tchi = {\rm Ker} ( H^d_m (I) \xrightarrow{\delta} H^d_Z(\mathcal{O}_X(\llcorner f^\ast D \lrcorner ))),$$

\noindent
where $\delta$ denotes an edge map of the Leray spectral sequence $H^i_m(H^j(X,\mathcal{O}_X(\llcorner f^\ast D \lrcorner))) \Rightarrow H^{i+j}_Z (\mathcal{O}_X( \llcorner f^\ast D \lrcorner)) $.

\end{Thm}

\begin{Rem}
{\rm The map $\delta$ is the degree one part of the graded map $\delta ' : H^d_n(S) \to H^d_{Z'}(\mathcal{O}_Y)$, which is also an edge map of the spectral sequence  $H^i_m(H^j(Y,\mathcal{O}_Y)) \Rightarrow H^{i+j}_{Z'} (\mathcal{O}_Y) $.The kernel of the map $\delta '$ equals $\tc _{H^d_m(S)}$. For an explanation of these claims, we refer the reader to Hara's paper,~\cite{ha}.}

\end{Rem}

Let us define now the multiplier ideal. 

\begin{Def}
\label{multiplier}
{\rm Let $V$ be an irreducible normal $\mathbf{Q}$-Gorenstein variety defined over a field $k$ of characteristic zero. Let $D=K_V+ \Delta$ be a $\mathbf{Q}$-divisor on $V$ and $X \xrightarrow{f} V$ be a proper birational map such that $f^*\Delta +E$ is a simple normal crossing divisor.} The multiplier ideal sheaf associated to $(V,D)$ {\rm is defined as

$$\mathcal{J}(V,\Delta) = f_*\mathcal(O_X(\ulcorner K_{X}- f^*D \urcorner))$$.

When $D$ is effective, the definition gives an ideal sheaf. In the general case, one has a fractional ideal sheaf. 

We would like to study the multiplier ideal in positive characteristic. 
To be able to define the multiplier ideal in characteristic $p \gg 0$, we start in characteristic zero and reduce the data to characteristic $p$. As mentioned above, this has been explained at length in ~\cite{ha}. To summarize the procedure, we indicate briefly its main points for the case $V={\rm Spec}(R)$, with $R$ finitely generated algebra over $k$ a field of characteristic zero (the case that we will be using later). We keep the notations just introduced in Definition~\ref{multiplier}. First choose a finitely generated $\mathbf{Z}$-subalgebra $A$ of $k$ and construct a finitely generated flat $A$-algebra $R_A$, a smooth $A$-scheme $X_A$, a birational morphism $f_A: X_A \xrightarrow {\rm Spec}(R_A)$, together with $\Delta_A$ $\mathbf{Q}$-divisor on ${\rm Spec}(R_A)$ and $E_A$, the exceptional fiber, such that $f_A^*\Delta_A +E_A$ has simple normal crossing and, by tensoring back with $k$, one obtains the initial data ${\rm Spec}(R), X, f, \Delta, E$. By choosing a general closed point $s$ in ${\rm Spec}(A)$ one gets the corresponding fibers $R_s,X_s,f_s$ etc., and the data are defined now over the residue field at $s$ which is of positive characteristic $p$. With all these data at hand one can define the multiplier ideal in characteristic $p \gg 0$ as above, in manner similar to the characteristic zero case.

For more details on multiplier ideals in characteristic zero, please consult ~\cite{la}.}
\end{Def}

The Theorem~\ref{hara} gives the following Corollary. We would like to recall that, for $M$ an Artinian $R$-module, $\tau_M$ denotes the test module earlier defined.

\begin{Cor}
Let $(R,m)$ be a local, normal, complete, $\mathbf{Q}$-Gorenstein of dimension d of characteristic $p \gg 0$ (obtained by reduction from characteristic zero). Using the same notations and hypotheses as above, 

$$\tau _{H^d_m(I) } = H^0(X, \mathcal{O}_X(\ulcorner K_{X}- f^\ast D \urcorner))$$

\noindent
seen as submodule of $\Hom_R(I, \omega _R)$ via the natural inclusion.

In particular, the multiplier ideal $H^0(X, \mathcal{O}_X(\ulcorner K_{X}- f^\ast D \urcorner))$ is a  strong test module.

\end{Cor}

\begin{Rem}
{\rm The proof of the first assertion of the Corollary follows closely the proof of Theorem 5.9 in Hara (which represents in fact the case $D=K_R$).}

\end{Rem}

\begin{proof}

By Hara's Theorem, 

$$\tchi = {\rm Ker}(H^d_m (I) \xrightarrow{\delta}  H^d_Z(\mathcal{O}_X(\llcorner f^\ast D \lrcorner ))).$$

Now, $$H^d_Z(\mathcal{O}_X(\llcorner f^\ast D \lrcorner )) = H^d_Z(\omega_X(- \ulcorner K_X - f^\ast D \urcorner )),$$ 

\noindent
which is Matlis dual to $H^0 (X, \mathcal{O}
_X(\ulcorner K_X - f^\ast D \urcorner ))$. Also, $H^d_m(I)^{\vee} = \Hom_R(I, \omega_R)$.

There is a natural inclusion of $H^0 (Y, \omega _Y) \into \omega _S$. Its degree one part gives a natural inclusion $H^0 (X, \mathcal{O}
_X(\ulcorner K_X - f^\ast D \urcorner )) \into \Hom_R(I, \omega_R)$.

We have that $$0 \to H^0 (X, \mathcal{O}
_X(\ulcorner K_X - f^\ast D \urcorner )) \to \Hom_R(I, \omega_R) \to \frac{\Hom_R(I, \omega_R)}{H^0 (X, \mathcal{O}
_X(\ulcorner K_X - f^\ast D \urcorner ))} = C \to 0.$$

By taking the Matlis dual, we get that
$$\tchi = C^\vee= {\rm Ann}_{H^d_m(I)} H^0 (X, \mathcal{O}
_X(\ulcorner K_X - f^\ast D \urcorner )).$$ Using again duality (as in Lemma~\ref{ann}), we get the first part of the Corollary. Theorem~\ref{strong} can be used now to conclude the proof.
\end{proof}

\begin{Rem}
{\rm  For each multiplier ideal as above, its minimum number of generators will provide a uniform bound (depending only on $R$) on the degree of the integral dependence equation of $x $ over $I$, for each $x \in I^\ast \subset \overline{I}$.}

\end{Rem}

\section{A few remarks on the parameter test ideal}

Many of the questions in tight closure theory that address the test ideal can also be formulated for an alternate notion of test ideal,  the parameter test ideal. Generally,  considering parameter ideals instead of arbitrary ideals provides questions with answers that have bearing on arbitrary ideals. In fact, tight closure is better understood in the case of ideals generated by parameters and many fundamental conjectures have been proven in these particular case (see for example Theorem 5.1 in~\cite{sm1}).

This final section deals with two natural questions that arise in the study of strong test ideals. The questions regard the parameter test ideal, so we will proceed by defining it (see Definition 8.7 in~\cite{ho-hu2}).

\begin{Def}
{\rm Let $R$ be an equidimensional local ring of positive
characteristic $p$. We define the} parameter test ideal $\taup(R)$
{\rm to be $\cap_{I} (I : I^\ast)$, where $I$ runs through all ideals
generated by a system of parameters.}
\end{Def}

Let us recall the basic properties of the parameter test ideal as in ~\cite{ho-hu2}.

\begin{Prop}
Let $(R,m,K)$ be an excellent equidimensional local ring of
characteristic $p$.

(a) If $c \in R^o$, then $c \in \taup(R)$ if and only if for every ideal
generated by a system of parameters, for all $x \in R$ we have that $x
\in I^\ast$ implies $cx^{[q]} \in I^{[q]}$ for all $q$.

(b) If $R$ is Cohen-Macaulay, $x_1, ..., x_d$ is an s.o.p. and $I_t =
(x_1^t,...,x_d^t)R$, then $\taup(R) = \cap_{t} (I_t : I_t^\ast)$.

(c) If $R$ is Cohen-Macaulay and $I$ is any ideal generated by
monomials in a system of parameters, then $\taup(R) (I^\ast) \subset
I$.

\end{Prop}

The parameter test ideal and its elements have also been studied by K. E. Smith (\cite{sm2}) and
J. V\'elez (\cite{ve}) with focus on Cohen-Macaulay excellent rings.  In particular,
K. E. Smith has shown that $\taup(R) = \Ann_{R} \hr = \Ann_{R} 0^{\ast
fg} _{H_{m}^{d}(R)}$ (Proposition 4.2,~\cite{sm2}).

The two questions are:

\noindent
\begin{Q}
\label{q1}
{\rm Is the parameter ideal
a strong test ideal for the family parameter ideals? More precisely, is it true that $\taup(R) I^\ast = \taup I$ for all ideals $I$ generated by systems of parameters? }
\end{Q}

In response to Question~\ref{q1}, we prove that $\taup(R) I^\ast = \taup I$ for a large family of ideals $I$ generated by systems of parameters.

\noindent
\begin{Q}
\label{q2}
{\rm  It is known that if $R$ is complete, then $\tc _{E_R(k)} = \Ann_{E_R(k)} (\tau(R))$. It is true that $\hr = \Ann_{H^d_m(R)} (\taup(R))$?}

\end{Q}

Denote by $N: = \Ann_{H^d_m(R)} (\taup(R))$. Clearly, $\hr \subset N$. (Whenever $R$ is excellent and analitically irreducible, $\hr$ is the unique maximal proper $F$-stable submodule of $H^d_m(R)$ as shown in~\cite{sm-thesis}). An affirmative answer to our Question~\ref{q2} would imply that $N$ is $F$-stable, because $\hr$ is $F$-stable. The problem of the $F$-stability of $N$ has appeared in the work of
K.~E.~Smith (~\cite{sm2}) and has remained open. ($N$ is $F$-stable for
Gorenstein rings, as the parameter test ideal equals the
test ideal in that case.)  If $N$ is $F$-stable, then $\taup(R)$ is what is called an {\it
$F$-ideal} of $R$. (For more on the notion of $F$-ideals, see~\cite{sm1,sm2}).  
We produce an example that settles the issue raised by Question~\ref{q2}, showing that  $\taup(R)$
is not necessarily an $F$-ideal, whenever $R$ is complete, Cohen-Macaulay and
reduced\footnote{After discussing the contents of our paper with Nobuo Hara, he informed us that, using some geometrical constructions, he can obtain an example that is not only reduced, but also normal, and for which $N$ is not $F$-stable.}.

We would like to show now that $\taup(R) I^\ast = \taup(R) I$ for {\it
certain} ideals generated by systems of parameters under some
additional conditions on the ring $R$. As previous authors did, we will concentrate on the case when $R$ is Cohen-Macaulay. First we need to state the
following:

\begin{Lma}Let $(R,m,K)$ be a local Cohen-Macaulay ring
and suppose that $\depth(\taup(R)) \geq 2$. Then for every ideal $I =
(c_1 d_1,c_2 d_2,x_3,...,x_d)$ generated by a system of parameters,
with $c_1,c_2$ are parameter test elements, $$I^\ast \subset I + (d_1
d_2) (c_1 ,c_2,x_3,...,x_d)^\ast.$$

\end{Lma}

\begin{proof}
Because each $c_i$ is a parameter test element and $I$ is generated by
a system of parameters, we have $c_i I^{\ast} \subset I$. So,
$I^{\ast} \subset I : (c_1,c_2)$. It can be shown that $I: (c_1,c_2)=
( d_1 d_2, x_3,...,x_d)$.

Each element $j$ of $I$ can be then  written in the form $j = i + d_1
d_2 r $, where $i \in I$ and $r \in R$. This gives that $ r \in I^\ast
: d_1 d_2$.   For every $q$ we have that for some $c$, a test element,
$c(r d_1 d_2)^{[q]} \in I^{[q]}$. Using the properties of regular
sequences, it follows that $cr^q \in (c_1, c_2, x_3,...,x_d)^{[q]}$,
for every $q \gg 0$. So, $r \in (c_1 ,c_2, x_3,...,x_d)^{\ast}$ and
this ends our proof.
\end{proof}

Now we can state one of the main results of this section:

\begin{Thm}
\label{strongpar}
If $(R,m,K)$ is Cohen-Macaulay and $\depth(\taup(R)) \geq 2$, then
there exist a large family of ideals $I$ generated by systems of parameters such that
$$\taup(R) I^\ast = \taup(R) I.$$

\noindent
More precisely, $\taup(R) I^\ast = \taup(R) I$, for every ideal parameter ideal  $I = (c_1 d_1,c_2 d_2,x_3,...,x_d)$, where $c_i, d_i$ are parameter test elements.

\end{Thm}

\begin{proof} Keep all the notations introduced in the above Lemma. Let us assume that $I = (c_1 d_1,c_2 d_2,x_3,...,x_d)$ is generated by a system of parameters with $c_1,c_2,d_1,d_2$ parameter test elements,

According to our Lemma, $$\taup(R) I^{\ast} \subset  \taup(R) I +
(d_1 d_2) \taup (R)(c_1,c_2,x_3,...,x_d)^{\ast}.$$  We have that
$$(d_1 d_2) \taup (R)(c_1,c_2,x_3,...,x_d)^{\ast} \subset (d_1
d_2)(c_1,c_2,x_3,...,x_d)$$ because the parameter test ideal has the
property that it multiplies the tight closure of every parameter ideal
into the ideal itself.

In conclusion, $$\taup(R) I^\ast \subset \taup(R) I+ (d_1
d_2)(c_1,c_2,x_3,...,x_d) \subset \taup(R) I.$$  using here in an
essential way that $d_1,d_2$ belong to $\taup(R)$.

So, $\taup(R) I^{\ast} \subset \taup(R) I.$

\end{proof}

The following Proposition is similar to results of Hara and Smith who proved that
the test ideal $\tau (R)$ is a strong test ideal whenever $\tau(R)=m$ (Theorem 1.1 in~\cite{ha-sm}) 

\begin{Prop}
\label{strongpar2}
Let $(R,m,k)$ Cohen-Macaulay and assume that $\taup(R) = m$. Then $m I^\ast = mI$ for all ideals $I$ generated by systems of parameters.
\end{Prop}

\begin{proof}
Fix $I$ an ideal generated by a system of parameters and assume that $m I^\ast \neq mI$. This means that there exist
$x \in I^\ast$ and $y \in m$ such that $yx \notin mI$. However, $\taup =m $ and $y \in m$, therefore $yx \in I$. Since $yx \in I-mI$, it follows that $yx$ can be taken as part of a minimal system of generators for $I$, say $I = (yx, x_2,...,x_d)$. 
Clearly, $yx, x_2,..., x_d$ form a system of parameters, hence they are a regular sequence in $R$.

Let us write that $x \in I^\ast$, by using $y$ as parameter test element. We get that, for every $q$, there exist $\lambda_q$ such that
$yx^q - \lambda_q y^qx^q \in (x_2,...,x_d)^{[q]}$. That is, $x^q \cdot ( y - \lambda_q y^q) \in (x_2,...,x_d)^{[q]}$.
But, $x, x_2,...,x_d$ form a regular sequence in $R$. So, $ y - \lambda_q y^q \in (x_2,...,x_d)^{[q]}$, and hence, 
$y \in (y,x_2,...,x_d)^{[q]}$ for every $q$. So, $y =0$ which is impossible.

\end{proof}

Theorem~\ref{strongpar} and Propostion~\ref{strongpar2} make us believe that Question~\ref{q1} has an affirmative answer for Cohen-Macaulay rings.

Now, we would like to study our second question: Is $\Ann_{H^d_m(R)} (\taup(R))
= \hr = \fghr$, where $(R,m,K)$ is a Cohen-Macaulay ring. (It is known
that $\hr = \fghr$, but we will not use this here.)

Let us denote by $N:= \Ann_{H^d_m(R)}(\taup(R))$. It is clear that
$\hr \subset N$. If one has equality, then it follows  that $N$ is
$F$-stable in $H^d_m(R)$, because $\hr$ is a proper submodule of
$H^d_m(R)$ that is stable under the action of Frobenius. We will
provide an example that shows that this is not true in
general. However, our example is not analytically irreducible (as we mentioned earlier, for rings that are excellent and analitically irreducible $\hr$ has the feature that it is a maximal proper $F$-stable submodule in $H^d_m(R)$).

Let $K$ be a field of characteristic $p$ and $x,y,z,w$ indeterminates
over $K$. Let $R$ be the ring $\frac{K[[x,y,z,w]]} {(xy,yz,zw)}$. This local
ring is the completion of the Stanley-Reisner ring
$K[x,y,z,w]/(xy,yz,zw)$ at the maximal ideal $(x,y,z,w)$. $R$ is
reduced, not domain with minimal primes $P_1=(x,z),\ P_2=(y,z),\
P_3=(y,w)$. Tight closure theory for Stanley-Reisner rings has been
studied by many authors and formulae for the parameter test ideal and
the test ideal have been given in terms of the minimal primes of the
ring (see for example Theorem 3.7 in~\cite{co}). We can apply this to $R$ and obtain $\taup(R) = \tau(R) = P_1
\cap P_2 + P_2 \cap P_3 + P_1 \cap P_3$. In our case, $P_1 \cap P_2 =
(xy,z)$, $ P_2 \cap P_3 = (y,zw)$ and $P_1 \cap P_3 = (xy,xw,zy,zw)$
and hence $\taup(R) = \tau(R) = (y,z,xw)$. Let us remark that $\Dim(R)
= 2$ and $\height (\taup(R)) = 1$.

Our ring $R$ is Cohen-Macaulay with a system of parameters given by
$x-w, \ x-y-z$. Denote by $I$ the ideal generated by these elements
and note that $R/I = K[y,z]/(y^2,yz,z^2)$.

We would like to show that there exists an element $\eta \in N$ for
which the image via the Frobenius action on $H^d_m(R)$, $F(\eta)$,
does not belong to $N$. This will show that $N$ is not $F$-stable and
therefore that it contains $\hr$ strictly.  Take $\eta = [(xw)^{p-1} +
(x^p-w^p, x^p -y^p-z^p)]$. We regard $H^d_m(R)$ as $\varinjlim_t
R/I_t$, where $I_t = (x^t-w^t,x^t-y^t-z^t)$ and the maps $R/I_t \to
R/I_{t+1}$, in our direct system, are given by multiplication by
$(x-w)(x-y-z)$. Now, $$\eta = [(xw)^{p-1} + (x^p-w^p, x^p -y^p-z^p)]$$
is written as an element of $R/(x^p-w^p, x^p -y^p-z^p) \subset
H^d_m(R)$, and $$F(\eta) = [(xw)^{p(p-1)} + (x^{p^2}-w^{p^2}, x^{p^2}
-y^{p^2}-z^{p^2})]$$ is written  as an element of $R/(x^{p^2}-w^{p^2},
x^{p^2} -y^{p^2}-z^{p^2}) \subset H^d_m(R)$.
 
To show that $\taup(R) \eta = 0$ we need to show that $xw \eta=0$,
because it is already clear that $y \eta = z \eta =0$, as $yx=zw = 0$
in $R$. It remains to check that $$(xw)^p \in (x^p-w^p, x^p -y^p-z^p,
xy, yz,zw)$$ in $K[[x,y,z,w]]$.  Since $yx \in (x^p-w^p, x^p -y^p-z^p,
xy, yz,zw)$, we get that $w^p y = x^py -y(x^p-w^p) \in (x^p-w^p, x^p
-y^p-z^p, xy, yz,zw)$. Since $w^py, wz$, and $w^p (x^p -y^p-z^p)  \in
(x^p-w^p, x^p -y^p-z^p, xy, yz,zw)$, we get that $(wx)^p \in (x^p-w^p,
x^p -y^p-z^p, xy, yz,zw)$ which is our claim.

Now we need to show that $F(\eta) \notin N$, that is $\taup(R) F(\eta)
\neq 0$. We will show that $$xw(xw)^{p(p-1)} \notin  (x^{p^2}-w^{p^2},
x^{p^2} -y^{p^2}-z^{p^2}, xy,yz,zw)$$ in $K[[x,y,z,w]]$. Let us assume
the contrary and take $y=z=0$.  This gives that $$xw(xw)^{p(p-1)} \in
(x^{p^2}-w^{p^2}, x^{p^2})=(x^{p^2},w^{p^2})$$ in $K[[x,w]]$ which is
certainly impossible.

This concludes our proof that $N$ is not $F$-stable in this example.

{\bf Acknowledgments:} I would like to thank Nobuo Hara, Mel Hochster and Paul Roberts for their helpful comments.

\end{document}